\documentclass[11pt,oneside]{article}
\usepackage[utf8]{inputenc}
\usepackage[english]{babel}
\usepackage{amsmath}
\usepackage{amsfonts}
\usepackage{amssymb}
\usepackage{amsthm}
\usepackage{tikz-cd} 
\usepackage{hyperref}
\usepackage{pgfplots}

\usepackage[left=3cm,right=3cm,top=3cm,bottom=3cm]{geometry}
\usepackage[title]{appendix}
\usepackage[]{algorithm2e}

\usepackage[normalem]{ulem}
\useunder{\uline}{\ul}{}
\theoremstyle{plain}
\newtheorem{teo}{}[section]
\newtheorem{prop}[teo]{Proposition}
\newtheorem{cor}[teo]{Corollary}

\newtheorem{lem}[teo]{Lemma}
\newtheorem{thm}[teo]{Theorem}
\theoremstyle{definition}
\newtheorem{ex}[teo]{Example}

\newtheorem{df}[teo]{Definition}

\newcommand\blfootnote[1]{%
  \begingroup
  \renewcommand\thefootnote{}\footnote{#1}%
  \addtocounter{footnote}{-1}%
  \endgroup
}

\definecolor{cof}{RGB}{219,144,71}
\definecolor{pur}{RGB}{186,146,162}
\definecolor{greeo}{RGB}{91,173,69}
\definecolor{greet}{RGB}{52,111,72}

\usepackage{graphicx}

\title{Semiflows on finite topological spaces}
\author{Pedro J. Chocano}
\date{}

\begin{document}
\maketitle

\begin{abstract}
In this paper, we study flows and semiflows defined on any given finite topological $T_0$-space $X$. We show that there exist non-trivial semiflows on $X$, unless $X$ is a minimal finite space. Specifically, non-trivial semiflows exist if and only if $X$ contains down beat points, and a non-trivial semiflow is essentially a strong deformation retraction. As a consequence of this result, we provide a new and concise proof that the only flow that can be defined on $X$ is the trivial flow. Finally, we discuss the number of different semiflows that can be defined on $X$ in terms of down beat points and other special points.

\end{abstract}
\blfootnote{2020  Mathematics  Subject  Classification:  	37C10, 54C15, 06A06 	 }
\blfootnote{Keywords: finite topological spaces, flows, semiflows, strong deformation retractions}
\blfootnote{This research is partially supported by Grant PID2021-126124NB-I00 from Ministerio de Ciencia, Innovación y Universidades (Spain).}
\section{Introduction}\label{sec:introduccion}

Recent research recognizes the importance of studying dynamical systems within the combinatorial context of finite partially ordered sets or finite topological $T_0$-spaces (see \cite{barmak2020Lefschetz,chocano2024coincidence} and the references therein). The idea of using these combinatorial spaces to approach problems in dynamical systems lies in the fact that they are computationally tractable, possess rich algebraic structure (see \cite{mccord1966singular}) and can be employed to reconstruct the phase space of dynamical systems (see \cite{chocano2020computational}). For these reasons, they are excellent candidates for developing computational methods to address problems in dynamical systems.

However, as pointed out in \cite[Proposition 2.4]{chocano2024coincidence}, the only flow that can be defined on a finite topological $T_0$-space is the trivial flow. Nevertheless, it is easy to find examples where there are non-trivial semiflows. In this paper, we focus on the study of semiflows defined on finite topological $T_0$-spaces. We prove that the existence of non-trivial semiflows depends on the homotopy type of the space $X$. Specifically, there exists a non-trivial semiflow on a finite topological space $X$ if and only if $X$ contains down beat points. These points have the particularity that, upon removal, the resulting space's homotopy type remains unchanged (see \cite{stong1966finite}). Moreover, we analyze the behaviour of non-trivial semiflows. Essentially, all of them are strong deformation retractions that move high points to lower points. Notice that there exists a general result stated in terms of finite sets and multivalued maps (see \cite[Theorem 3.1]{barmak2024conley}) which also indicates that a semiflow is equal for positive values, but it does not provide specific information about the semiflow on the positive values.

The organization of the paper is as follows. In Section \ref{Section:preliminaries}, we establish the basic concepts from the theory of finite topological spaces needed to study semiflows within this context. In Section \ref{section:semiflows}, we illustrate the existence of non-trivial semiflows through a concrete example. Although  non-trivial semiflows exist, we prove that an important family of finite topological spaces does not admit non-trivial semiflows. From this and other results, we conclude that the existence of non-trivial semiflows is related to the homotopy type of the spaces on which the semiflow is defined. To conclude this section, we apply the results obtained to provide a concise proof that the only flow on a finite topological $T_0$-space is the trivial flow. Finally, in Section \ref{sec_counting_semiflows}, we discuss and obtain results regarding to the number of distinct semiflows that can be defined on a finite topological space $X$ in terms of down beat points and other special points.

\section{Preliminaries}\label{Section:preliminaries}

Let $X$ be a finite topological $T_0$-space, and let $x\in X$. We denote by $U_x$ ($F_x$) the intersection of every open (closed) set containing $x$. It is evident that $U_x$ ($F_x$) is an open (closed) set. Given two points $x,y\in X$, we say that $x\leq y$ if and only if $U_x\subseteq U_y$. The set $X$ with this relation forms a finite partially ordered set (poset for short). Conversely, given a finite partially ordered set $(X,\leq)$, the set of lower sets forms a basis for a $T_0$ topology on $X$. Recall that $S\subseteq X$ is a lower set if whenever $y\leq x$ and $x\in S$, then $y\in S$. These two association are mutually inverse, demonstrating that for any finite set $X$, there exists a bijective correspondence between the $T_0$ topologies and the partial orders defined on $X$. From now on, we treat finite posets and finite topological $T_0$-spaces as equivalent objects without explicit mention. Additionally, we assume that every finite topological space is a $T_0$-space. It is worth noting that a map $f:X\rightarrow Y$ between finite topological spaces is continuous if and only if it is order preserving. From this, it can be deduced that the category of finite posets and the category of finite topological spaces are isomorphic \cite{alexandroff1937diskrete}. We will assume, without explicit mention, that every finite topological space that appears throughout this paper is a $T_0$-space.

The \emph{height} $ht(X)$ of a finite poset $X$ is one less than the maximum number of elements in a chain of $X$. The \emph{height of a point} $x\in X$ is $ht(x)=ht(U_x)$. Given a finite partially ordered set $(X,\leq)$, and $x,y\in X$, we write $x\prec y$ when $x<y$ and there is no element $z\in X$ such that $x<z<y$. The \emph{Hasse diagram} of $X$ is a directed graph whose vertices are the points of $X$, and there is an edge from $x$ to $y$ if and only if $x\prec y$. In subsequent Hasse diagrams we omit the orientation of the edges. Let $f,g:X\rightarrow Y$ be continuous maps between finite topological spaces, we write $f\leq g$ whenever $f(x)\leq g(x)$ for every $x\in X$. We now recall and introduce basic concepts from the theory of finite topological spaces. For a comprehensive introduction, we refer \cite{barmak2011algebraic}.

\begin{prop}\cite[Corollary 1.2.6]{barmak2011algebraic}\label{prop:homotopia} Let $f,g:X\rightarrow Y$ be two continuous maps between finite topological spaces. Then $f$ is homotopic to $g$ if and only if there is a fence $f=f_0\leq f_1\geq f_2\leq \cdots\geq f_n=g$.
\end{prop}

\begin{df} A point $x$ in a finite topological space $X$ is a down (up) beat point if $U_x\setminus \{x \}$ ($F_x\setminus \{x \}$) has a maximum (minimum).
\end{df}

\begin{df} Let $X$ be a topological space and $A\subseteq X$. A \emph{retraction} $r:X\rightarrow A$ is a continuous map satisfying that $r(a)=a$ for every $a\in A$. We say that $r$ is a \emph{deformation retraction} if it is a retraction and $i\circ r$ is homotopic to the identity map $\text{id}:X\rightarrow X$, where $i:A\rightarrow X$ denotes the inclusion map. The space $A$ is called \emph{deformation retract} of $X$. Additionally, if in the previous definition, the homotopy between $i\circ r$ and the identity map leaves fixed the points of $A$, then we say that $A$ is a \emph{strong deformation retract} of $X$. 
\end{df}

\begin{prop}\cite[Proposition 1.3.4]{barmak2011algebraic}\label{prop:homotopia_finitos} Let $X$ be a finite topological space and let $x\in X$ be a beat point. Then $X\setminus \{ x\}$ is a strong deformation retract of $X$.
\end{prop}

A finite topological space $X$ is a \emph{minimal finite space} if it does not have beat points.

\begin{ex} Consider $X=\{A,B,C,D,E\}$ with $A,B<C,D$ and $A,B,C<E$. The finite topological space $X$ is not a minimal finite space because $E$ is a down beat point. However, $X\setminus \{ E\}$ is a minimal finite space because it does not contain beat points. 
\end{ex}

\begin{thm}\cite[Theorem 1.3.6]{barmak2011algebraic}\label{thm:minimal_homotopica_identidad_es_identidad} Let $X$ be a minimal finite space. A continuous map $f:X\rightarrow X$ is homotopic to the identity map $\textnormal{id}:X\rightarrow X$ if and only if $f=\textnormal{id}$.
\end{thm}

To conclude, we recall the basic concept of this paper. 

\begin{df}\label{def:dynamical_system} By a \emph{semiflow}  on a topological space $X$ we mean a continuous map $\varphi:\mathbb{R}_0^+\times X\rightarrow X$ such that:
\begin{enumerate}
\item[i)] For every $x\in X$, $\varphi(0,x)=x$.
\item[ii)] For every $x\in X$ and $s,t\in \mathbb{R}_0^+$, $\varphi(s,\varphi(t,x))=\varphi(s+t,x)$.
\end{enumerate}
If the group $\mathbb{R}$ replaces the semigroup $\mathbb{R}_0^+$ of the non-negative real numbers, then $\varphi$ is called  \emph{flow}.
\end{df}

\textbf{Notation.} Let $\varphi$ be a semiflow (or a flow) and let $t$ be a non-negative real number (or a real number). We will denote by $\varphi_t$ the continuous map $\varphi(t,\cdot):X\rightarrow X$. We say that a semiflow (or a flow) is trivial if $\varphi_t$ is the identity map $\text{id}:X\rightarrow X$ for every $t\in \mathbb{R}_0^+$ (or $t\in \mathbb{R}$).

In \cite[Proposition 2.4]{chocano2024coincidence}, the authors proved that the concept of flow within the context of finite topological spaces or finite posets is not particularly interesting. Furthermore, in \cite{chocano2022dynamical}, they generalized this result to Alexandroff spaces or posets. For this reason, to study dynamical systems on finite topological spaces, the authors introduced multivalued maps. As an application of their study, they approximated classical dynamical systems using finite spaces and multivalued maps (see \cite[Section 6]{chocano2024coincidence}). .


\section{Semiflows on finite topological spaces}\label{section:semiflows}

We begin by proving the existence of non-trivial semiflows through an illustrative example. This example will subsequently highlight properties relevant to semiflows within this context.
\begin{ex}\label{Ex:notrivial} Let $X:=\{A,B,C,D,E,F \}$ with $A>B,C,D,E,F$; $B>D,E,F$; $C>D,E,F$, and $D>E,F$. Define $\varphi:\mathbb{R}_0^+\times X\rightarrow X$ by 

\begin{align*}
\varphi(t,x)&:=x \ \ \ \textnormal{if } t=0,\\
\varphi(t,x)&:= \begin{cases} D & \ \ \ \textnormal{if }x> D \  \textnormal{and }t\in (0,\infty), \\
x & \ \ \ \textnormal{if }x\leq  D \  \textnormal{and }t\in (0,\infty).
\end{cases}
\end{align*}
It is an easy exercise to prove that $\varphi$ is a semiflow. Moreover, $\varphi$ is a non-trivial semiflow because $\varphi_t\neq \text{id}$ for every $t\in \mathbb{R}^+$. We present in Figure \ref{fig:no_trivial} the Hasse diagram of $X$, along with a schematic representation of $\varphi$.
\begin{figure}[h!!!]
\centering
\includegraphics[scale=1]{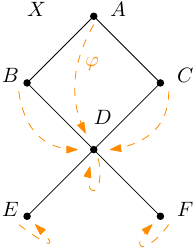}
\caption{Hasse diagram of $X$ and schematic representation of the semiflow $\varphi$ from Example \ref{Ex:notrivial}.}\label{fig:no_trivial}
\end{figure}

\end{ex}

\begin{thm}\label{thm:trivial_minimal} Let $X$ be a minimal finite topological space. If $\varphi:\mathbb{R}_0^+\times X \rightarrow X$ is a semiflow, then $\varphi$ is trivial.
\end{thm}
\begin{proof}
We argue by contradiction. Suppose that there exists a non-trivial semiflow $\varphi:\mathbb{R}_0^+\times X \rightarrow X$. This implies that there exist $x\in X$ and $t\in \mathbb{R}^+$ such that $\varphi(t,x)\neq x$. Moreover, the continuous map $\varphi_t:X\rightarrow X$ is homotopic to the identity map $\text{id}:X\rightarrow X$ because $\varphi_{|[0,t]}:[0,t]\times X\rightarrow X$ provides a homotopy between $\text{id}$ and $\varphi_t$. However, by Theorem \ref{thm:minimal_homotopica_identidad_es_identidad}, we obtain that $\text{id}=\varphi_t$, which entails the contradiction.
\end{proof}

Theorem \ref{thm:trivial_minimal} proves that we only have non-trivial dynamical systems within this framework unless the space $X$ contains beat points. For instance, in Example \ref{Ex:notrivial}, the points $B,C, E$ and $F$ are beat points, allowing us to define a non-trivial semiflow.  In the following sequence of lemmas, we state and narrow down the possibilities for semiflows.

\begin{lem}\label{lem:trivial_extiende_identidad_fuera_cero} Let $X$ be a finite topological space and let $\varphi:\mathbb{R}_0^+\times X\rightarrow X$ be a semiflow. If there exists $\epsilon>0$ such that $\varphi_t$ is the identity map for every $t\in [0,\epsilon)$, then $\varphi$ is trivial.
\end{lem}
\begin{proof}
Let $s\in \mathbb{R}_0^+$. Then $s=\sum_{i=1}^n {s_i}$, where $s_i\in (0,\epsilon)$ for every $i=1,...,n$. Consequently, by ii) in Definition \ref{def:dynamical_system}, for every $x\in X$, we obtain
$$\varphi_s(x)=\varphi_{s_n}\circ \cdots\circ \varphi_{s_2}\circ \varphi_{s_1}(x)=x, $$
and this shows that $\varphi$ is trivial.
\end{proof}

Observe that this lemma is in concordance with the definition of $\varphi$ from Example \ref{Ex:notrivial}.

\begin{lem}\label{lem:cadenas} Let $X$ be a finite topological space and let $\varphi:\mathbb{R}_0^+\times X\rightarrow X$ be a semiflow. If $s,t\in \mathbb{R}_0^+$ are such that $s<t$, then $\varphi_t\leq \varphi_s$. 
\end{lem}
\begin{proof}
The map $\varphi$ is continuous at $(0,x)$ and $(s,x)$ for every $x\in X$, which gives that there exist $\epsilon>0$ and $\delta>0$  satisfying that $\varphi([0,\epsilon)\times U_x)\subseteq U_{x}$ and $\varphi(I\times U_x)\subseteq U_{\varphi_s(x)}$ where $I=\mathbb{R}_0^+\cap (s-\delta,s+\delta)$. We can express $t$ as follows $t=\sum_{i=1}^n s_i+s_{n+1}$, where $s_i\in (0,\epsilon)$ and $s_{n+1}\in I$ for every $i=1,...,n$. By ii) in Definition \ref{def:dynamical_system}, we obtain
$$\varphi_t(x)=\varphi_{s_{n+1}}\circ \varphi_{s_n} \circ  \cdots \circ \varphi_{s_2}\circ \varphi_{s_1}(x).$$
It is evident that $y=\varphi_{s_n} \circ  \cdots \circ \varphi_{s_2}\circ \varphi_{s_1}(x)\in U_x$. Consequently, $\varphi_{t}(x)=\varphi_{s_{n+1}}(y)\in U_{\varphi_{s}(x)}$.
\end{proof}
As an immediate consequence of this result, we can obtain the following lemma.
\begin{lem}\label{lem:varios resultados} Let $X$ be a finite topological space and let $\varphi:\mathbb{R}_0^+\times X\rightarrow X$ be a semiflow. Then
\begin{enumerate}
\item $\textnormal{id}\geq  \varphi_t$ for every $t\in \mathbb{R}_0^+$.
\item $\varphi(t,x)\in U_x$ for every $t\in \mathbb{R}_0^+$ and $x\in X$.
\item $\varphi(t,y)=y$ for all $t\in \mathbb{R}_0^+$ and every  $y\in X$ with $ht(y)=0$.
\end{enumerate}
 
\end{lem}


These results show that the dynamics of a semiflow always move from top elements to lower elements. In the Hasse diagram of a finite topological space, we progress toward the bottom (see Figure \ref{fig:no_trivial}). Moreover, there must be down beat points to define non-trivial semiflows. 
\begin{thm} Let $X$ be a finite topological space. If $X$ does not have down beat points and $\varphi:\mathbb{R}_0^+\times X\rightarrow X$ is a semiflow, then $\varphi$ is trivial.
\end{thm}
\begin{proof}
Let $x\in X$. We will show that $\varphi_t(x)=x$ for every $t\in \mathbb{R}_0^+$ by induction on the height of $x$. Suppose $ht(x)=0$; then, by 3. in Lemma \ref{lem:varios resultados}, we have  $\varphi_t(x)=x$ for every $t\in \mathbb{R}_0^+$. The case $ht(x)=1$ follows trivially because $U_x\setminus \{x \}$ is a disjoint set of points, all of which are fixed by $\varphi_t$. Consequently, due to the continuity of $\varphi_t$ for every $t\in \mathbb{R}_0^+$, we conclude that $\varphi_t(x)=x$. Now let us assume that the result holds when $ht(x)=n$. Consider the case where $ht(x)=n+1$. By the hypothesis, $U_x\setminus \{x\}$ does no have a maximum, so there exist finitely many points $y_i$ such that $y_i\prec x$, where $i=1,...,m$. It is evident that $ht(y_i)<n+1$ for every $i=1,...,m$. By the induction hypothesis, we know that $\varphi_t(y_i)=y_i$ for every $t\in \mathbb{R}_0^+$ and $i=1,...,m$.  If $\varphi_t(x)\leq y_i<x$ for some $i=1,...,m$, then, due to the continuity of $\varphi_t$, we have $y_j=\varphi_t(y_j)\leq \varphi_t(x)\leq y_i$ for every $j\neq i$, which leads to a contradiction. Thus, $\varphi_t(x)=x$.
\end{proof}

Given a finite topological space $X$, let $D(X)$ denote the set of down beat points of $X$. Then, there are non-trivial semiflows on $X$ if and only if $D(X)$ is non-empty. For each $x\in D(X)$, we denote by $\overline{x}$ the maximum of $U_x\setminus \{ x\}$.

\begin{ex}\label{ex:semiflow_natural} Let $X$ be a finite topological space such that $D(X)$ is non-empty. Define $\varphi:\mathbb{R}_0^+\times X\rightarrow X$  by 
\begin{align*}
\varphi(t,x)&:=x \ \ \ \textnormal{if } t=0,\\
\varphi(t,x)&:= \begin{cases} \overline{x} & \ \ \ \textnormal{if }x\in D(X) \  \textnormal{and }t\in (0,\infty), \\
x & \ \ \ \textnormal{if }x\notin  D(X) \  \textnormal{and }t\in (0,\infty).
\end{cases}
\end{align*}
It is straightforward to verify that $\varphi$ is a non-trivial semiflow. Moreover, it is evident that $\varphi_t(X)$ is a strong deformation retract of $X$ for every $t\in (0,\infty)$. We deduce this result from Proposition \ref{prop:homotopia_finitos}. We can interpret $\varphi$ as a process of removing some of the down beat points of $X$. 
\end{ex}

Note that not every semiflow on a finite topological space $X$ is the same as the one defined in Example \ref{ex:semiflow_natural}. For instance, consider Example \ref{Ex:notrivial}. Furthermore, it is possible to have contractible spaces that do not admit non-trivial semiflows. To illustrate this, let $X$ be a non-contractible space. Define $Y$ as the union of $X$ and a single point $ \{y \}$ keeping the given ordering in $X$ and declaring that $y>x$ for every $x\in X$. Clearly, $Y$ is contractible, but $D(Y)=\emptyset$. 

\begin{thm}\label{thm_retracto_deformacion_fuerte} Let $X$ be a finite topological space and let $\varphi:\mathbb{R}_0^+\times X\rightarrow X$ be a non-trivial  semiflow. Then $\varphi_t=r$ for every $t\in \mathbb{R}^+$ where $r:X\rightarrow r(X)$ is retraction and $r(X)$ is a strong deformation retract of $X$.
\end{thm}
\begin{proof}
First, observe that there is a finite number of continuous maps $r:X\rightarrow X$ such that $r\leq \text{id}$, where $\text{id}$ is the identity map, because $X$ is a finite topological space. Thus, there are a finite number of chains of continuous maps of the following form $\text{id}>r_1>\cdots>r_n$. Note also that these chains are all finite. From this, we can deduce that there exists a positive real value $a$ satisfying that $\varphi_t=r:X\rightarrow X$ for every $t\in (0,a)$, where $r$ is a continuous map such that $\text{id}>r$. Otherwise, we would get a contradiction with the previous comments.

Now, we prove that $r(X)$ is a strong deformation retract of $X$ by proving that $i\circ r=\text{id}$ and that $r\circ i$ is homotopic to $\text{id}$ where $i:r(X)\rightarrow X$ denotes the inclusion. Let us consider $t\in (0,a)$, it is evident that $t=t_1+t_2$ for some $t_1,t_2\in (0,t)$. Then 
\begin{align*}
\varphi_t(x)=\varphi_{t_1}\circ \varphi_{t_2}(x),
\end{align*}
which is the same as $r(x)=r^{2}(x)$. This proves that $r$ is the identity map on $r(X)$, i.e., $r:X\rightarrow r(X)$ is a retraction. By the construction, we have $i\circ r\leq \text{id}$. This yields that  $r\circ i$ is homotopic to $\text{id}$ from Proposition \ref{prop:homotopia}. Thus, $r(X)$ is a strong deformation retract.

We will show that $\varphi_t=r$ for every $t\in \mathbb{R}^+$. Suppose $t\in [a,\infty)$. Then $t=\sum_{i=1}^n t_i$, where $t_i\in (0,a)$ for every $i=1,...,n$. Consequently, 
\begin{align*}
\varphi_t(x)& =\varphi_{t_n}\circ \cdots\circ \varphi_{t_1}(x) \\
& =r\circ \cdots\circ r(x)=r^n(x)=r(x).
\end{align*}
\end{proof}

As an immediate consequence of this result, we obtain a new proof of \cite[Proposition 2.4]{chocano2024coincidence}.

\begin{cor}\label{thm_no_flows} Let $X$ be a finite topological space. If $\varphi:\mathbb{R}\times X\rightarrow X$ is a flow, then $\varphi$ is trivial.
\end{cor}
\begin{proof}
We argue by contradiction. Let us suppose that $\varphi$ is non-trivial. Then there exists $t>0$ and $x\in X$ such that $\varphi(t,x)=y\neq x$ (if $t<0$, then $\varphi(t,x)= y\neq x$ and we have that $\varphi(-t,y)=x$, where $-t>0$). On the other hand, $\varphi_{|\mathbb{R}_0^+\times X}$ is a semiflow. By Theorem \ref{thm_retracto_deformacion_fuerte}, $\varphi(t,x)=r(x)$ for every $t>0$, where $r:X\rightarrow X$ is a strong deformation retraction. However, $r$ does not have inverse unless $r$ is the identity, which entails the contradiction.
\end{proof}

\begin{figure}[h!!]
\centering
\includegraphics[scale=1]{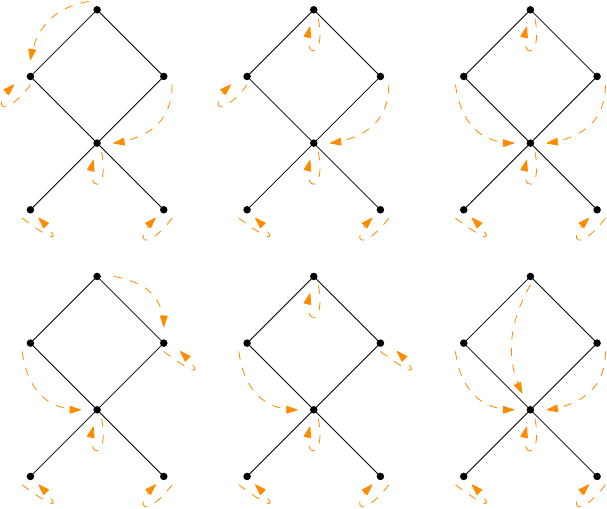}
\caption{Description of the non-trivial semiflows defined on the space $X$ from Example \ref{Ex:notrivial}.}\label{fig:todos_semiflujos}
\end{figure}

In conclusion, given a finite topological space $X$, a non-trivial semiflow $\varphi:\mathbb{R}_0^+\times X\rightarrow X$ is essentially a strong deformation retraction. Moreover, since $X$ is finite, there are only a finite number of semiflows that can be defined on it. For instance, in Example \ref{Ex:notrivial}, it is easy to verify that there are only 6 non-trivial semiflows that can be defined on $X$. We have represented them in Figure \ref{fig:todos_semiflujos}.

\section{Counting semiflows}\label{sec_counting_semiflows}

Given any finite topological space $X$, we define $S_F(X)$ as the number of distinct semiflows that can be defined on $X$. As a direct consequence of the results from Section \ref{section:semiflows}, we may derive the following theorem.

\begin{thm} Let $X$ be a finite topological space. Then $D(X)=\emptyset$ if and only if $S_F(X)=1$.  
\end{thm}

We now study relations between the number of down beat points of a finite topological space $X$ and the number $S_F(X)$.

\begin{prop}\label{prop_cota_down_beat_point} Let $X$ be a finite topological space. Then $S_F(X)\geq 2^{|D(X)|}$.
\end{prop}
\begin{proof}
Set $D(x)=\{x_1,...,x_n\}$ such that $ht(x_i)\geq ht(x_j)$ for every $1\leq i<j\leq n$. For each down beat point $x_i\in D(X)$, we denote by $\overline{x}_i$ the maximum of $U_{x_i}\setminus \{x_i\}$. Consider $r_i(x_i):=\overline{x}_i$ and $r_{i_{|X\setminus \{ x_i\}}}:=\textnormal{id}:X\setminus \{ x_i\}\rightarrow X\setminus \{ x_i\}$. We can obtain $2^n$ distinct semiflows by combining the retractions $r_i$ defined above. We prove the last assertion. Let $I=(a_1,...,a_n)\in \{0,1 \}^n$ and set $X^I_i:=X\setminus \{ x_j\in D(X)|j\leq i$ and $a_j=1\}$. For each $1\leq i \leq n$, define $r_{a_i}:X^I_{i-1}\rightarrow X^I_i$ by $r_{a_i}:=\textnormal{id}$ if $a_i=0$ and $r_{a_i}(x):=r_{i_{|X^I_{i-1}}}(x)$ if $a_i=1$, where $X_0^I=X$. Consider $r_I:X\rightarrow X^I_n$ defined by $r_I(x):=r_{a_n}\circ r_{a_{n-1}}\circ \cdots \circ r_{a_2}\circ r_{a_1}(x)$. We verify that $r_I$ is a strong deformation retraction. By the construction, it is clear that $r_{I}\circ i=\textnormal{id}_{X^I_n}$, where $i:X^I_n\rightarrow X$ is the inclusion, and also that $i\circ r_I\leq \textnormal{id}_X$. Given $I,J\in \{0,1\}^n$, it is evident that $r_I=r_J$ if and only if $I=J$. Hence, there are $2^n$ distinct strong deformation retractions of the form $r_I$. To conclude, we have that $\varphi^I :\mathbb{R}_0^+\times X\rightarrow X$, defined by $\varphi^I(0,x):=x$ with $x\in X$ and $\varphi^I(t,x):=r_I(x)$ for every $t>0$ and $x\in X$, is a semiflow.
\end{proof}

\begin{figure}[h!]
\centering
\includegraphics[scale=1]{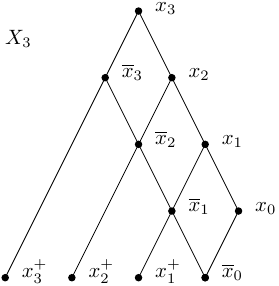}
\caption{Hasse diagram of $X_3$ from proof of Theorem \ref{thm_realizacion_todos_los_enteros_positivos}.}\label{fig:todos_numeros_realizados}
\end{figure}

\begin{thm}\label{thm_realizacion_todos_los_enteros_positivos} For every positive integer number $n\geq 2$ there exists a finite topological space $X$ such that $D(X)=1$ and $S_F(X)=n$. 
\end{thm}
\begin{proof}
Consider the finite topological space $X_n:=\{\overline{x}_0, x_0, x_1^+,\overline{x}_1,x_1,x_2^+,\overline{x}_2,x_2,...,x_n^+,\overline{x}_n,x_n \}$ given by $\overline{x}_i<x_{j},\overline{x}_j$; $x_i<x_j$ and $x_i^+<x_j,\overline{x}_j$ for every $i<j$, where $n$ is a non-negative integer number (see the Hasse diagram of $X_3$ in Figure \ref{fig:todos_numeros_realizados}). It is evident that $X_n$  has only one down beat point, $D(X_n)=\{ x_0\}$. Moreover, it is simple to deduce that $\{x_0,x_1,...,x_n \}$ are the only points that are down beat points of $X$ after removing down beat points. Consider the retractions $r_i:X_n\rightarrow X_n$ defined by $r_i(x_j):=\overline{x}_j$ for all $0\leq j\leq i$ and $r_i(x):=x$ otherwise, where $i\leq n$. The non-trivial semiflows on $X$ are given by $r_i$. Thus, $S_F(X_{n})=n+2$. 
\end{proof}

Note that by modifying the space $X_n$ from the proof of Theorem \ref{thm_realizacion_todos_los_enteros_positivos} we can produce infinitely many finite topological spaces with $n$ distinct semiflows. The construction proposed in Theorem  \ref{thm_realizacion_todos_los_enteros_positivos} shows that down beat points are not the unique important points to obtain non-trivial semiflows. In this case, the points $x_i$ from the space $X_3$ of Figure \ref{fig:todos_numeros_realizados} play an important role. In the following definition we identify that sort of points.

\begin{df}\label{def_potential_down_beat_point} Let $X$ be a finite topological space and $y\in X$. We say that $y$ is a \emph{potential down beat} point if there exists a sequence of points $x_1,...,x_n\in X$ such that $x_1$ is a down beat point of $X$, $x_{i+1}$ is a down beat point of $X\setminus \{x_1,...,x_i\}$ for every $i=1,...,n-1$, $ht(x_i)\leq ht(x_j)$ in $X$ if and only if $i\leq j$ and $y=x_n$.
\end{df}

Given a finite topological space $X$, it is evident that every down beat point $x\in X$ is a potential down beat point by considering the trivial sequence: $x$. Moreover, if $x$ is a potential down beat point, then there exists a down beat point $y\in X$ such that $y<x$.

\begin{thm}\label{thm:los_flujos_mueven_potenciales_beat_points} Let $X$ be a finite topological space and $y\in X$. Then, there exists a non-trivial semiflow $\varphi:\mathbb{R}_0^+\times X\rightarrow X$  such that $\varphi(t,y)\neq y$ for every $t>0$ if and only if $y$ is a potential down beat point.
\end{thm}
\begin{proof}

Suppose that $y$ is a potential down beat point. Then there exists a sequence $x_1,...,x_n\in X$ such that $x_1$ is a down beat point of $X$, $x_{i+1}$ is a down beat point of $X\setminus \{x_1,...,x_i\}$ for every $i=1,...,n-1$, $ht(x_i)\leq ht(x_j)$ in $X$ if and only if $i\leq j$, and $y=x_n$. Let us denote by $\overline{x}_{i}$ the maximum of $U_{x_i}$ in $X\setminus \{x_1,...,x_{i-1} \}$ for $i>2$ and let $\overline{x}_{1}$ be the maximum of $U_{x_1}$ in $X$.  Define $\varphi:\mathbb{R}_0^+\times X\rightarrow X$ by $\varphi(0,x):=x$ and $\varphi(t,x):=r(x)$, where $r:X\rightarrow X$ is the retraction given by $r(x_{i}):=\overline{x}_i$ for every $i=1,...,n$ and $r_{|X\setminus \{x_1,...,x_{n}\}}:=\text{id}:X\setminus \{x_1,...,x_{n}\}\rightarrow X\setminus \{x_1,...,x_{n}\}$.  By the construction of $r$, it is clear that $r(x)\leq x$ for every $x\in X$. To prove the continuity of $r$, we study the possible cases that we have where we consider $x<x'$ in $X$. Suppose $x,x'\notin \{ x_1,...,x_n\}$. Then $x=r(x)<r(x')=x'$. Suppose that $x=x_i$ for some $i=1,...,n$ and $x'\notin \{ x_1,...,x_n\}$. This implies that $r(x_i)=\overline{x}_i<x_i<x'=r(x')$. Suppose $x\notin \{ x_1,...,x_n\}$ and $x'=x_i$ for some $i=1,...,n$. By the definition of $r$, $r(x)=x$ and $r(x_i)=\overline{x}_i$. If $\overline{x}_i\not \sim x$, then we get that $x_i$ is not a down beat point in $X\setminus \{x_1,...,x_{i-1}\}$, and this implies a contradiction. If $x>\overline{x}_i$, then $\overline{x}_i$ is not the maximum of $U_{x_i}\setminus \{x_i \}$ in $X\setminus \{x_1,...,x_{i-1} \}$ because $\overline{x}_i<x<x_i$, which leads to a contradiction. Suppose $x,x'\in \{x_1,...,x_n \}$, let us say $x=x_i$ and $x'=x_j$. By Definition \ref{def_potential_down_beat_point}, we get that $i<j$ because $ht(x_i)<ht(x_j)$. By the construction of $r$, $\overline{x}_i=r(x_i)<x_i$ and $\overline{x}_j=r(x_j)<x_j$. Therefore, $\overline{x}_i<x_j$. If $\overline{x}_i\notin \{ x_1,...,x_n\}$, then $\overline{x}_i\leq \overline{x}_j$ since $\overline{x}_j$ is the maximum of $U_{x_j}$ in $X\setminus \{ x_1,...,x_{j-1}\}$. If $\overline{x}_i\in \{x_1,...,x_n \}$, then we get a contradiction with Definition \ref{def_potential_down_beat_point}. To conclude, we show that $r\circ i=\textnormal{id}_{r(X)}$ where $i:r(X)\rightarrow X$ is the inclusion. We prove that $r(x_i)\notin \{ x_i,...,x_n\}$ for every $i$. We argue by contradiction, suppose $r(x_i)=x_j$ for some $j>i$, i.e., $\overline{x}_i=x_j$. By the construction of $r$, we obtain $x_i>x_j$.  By the hypothesis, $i>j$ because $ht(x_i)> ht(x_j)$, and this leads to a contradiction. Moreover, by the construction of $r$, it is evident that $r(x_i)\neq x_j$ with $j\leq i$. Hence, $r(x_i)\in X\setminus \{x_1,...,x_n \}$ for every $1\leq i\leq n$. Since $r_{|r(X)}=\textnormal{id}_{r(X)}$, we obtain $r\circ i=\textnormal{id}_{r(X)}$.

Suppose that there exists a non-trivial semiflow $\varphi$ such that $\varphi(t,y)\neq y$. Define $r:X\rightarrow X$ by $r(x):=\varphi(1,x)$. Let us consider $\{x_1,...,x_n \}:=X\setminus r(X)$, where we have $i<j$ whenever $ht(x_i)\leq ht(x_j)$. By 1. in Lemma \ref{lem:varios resultados} and the fact that  $r(X)$ is a strong deformation retract of $X$ (see Theorem \ref{thm_retracto_deformacion_fuerte}), we may deduce that $y$ is a potential down beat point by considering  $\{x_1,...,x_n \}$.
\end{proof}



\begin{prop}\label{prop_cota_inferior_potential_down_beat_point_anticadena} Let $X$ be a finite topological space and let $A$ be the maximal antichain of potential down beat points of $X$ such that $U_a\cap U_b=\emptyset$ for every $a,b\in A$ with $a\neq b$. Then $S_F(X)\geq 2^{|A|}$.
\end{prop}
\begin{proof}
By the proof of Theorem \ref{thm:los_flujos_mueven_potenciales_beat_points} and the hypothesis, for each potential down beat point $a$ of $A$ there exists a non-trivial semiflow $\varphi^a$ such that $\varphi^a(t,a)\neq a$ for every $t>0$ and $\varphi^a(t,x)=x$ for every $t\geq 0 $ and $x\in X\setminus \{ U_a\}$. For every $a\in A$, set $r_a:=\varphi^a(1,x):X\rightarrow X$. By repeating the same arguments from the proof of Proposition \ref{prop_cota_down_beat_point} for the retractions $r_a$, we can obtain the desired result.
\end{proof}

\begin{prop}\label{prop_si_potential_se_mueve_un_down_tambien} Let $X$ be a finite topological space and $x\in X $. If $x$ is not a down beat point and $\varphi$ is a non-trivial semiflow defined on $X$ such that $\varphi(t,x)\neq x$ for every $t>0$, then there exists a down beat point $y<x$ such that $\varphi(y,t)\neq y$ for every $t>0$.
\end{prop}
\begin{proof}
Let $r:X\rightarrow X$ be the retraction that defines $\varphi$. There exist $y_1,y_1'\in X$ such that $y_1\prec x$ and $y_1'\prec x$ with $y_1\not \sim y_1'$, and $r(x)\leq y_1$ because $x$ is not a down beat point. By the continuity of $r$, we deduce that $r(y_1')\neq y_1'$. By Theorem \ref{thm:los_flujos_mueven_potenciales_beat_points}, $y_1'$ is a potential down beat point.  If $y_1'$ is a down beat point, then we can conclude. Suppose that $y_1'$ is not a down beat point. Then we can repeat the previous argument to obtain a potential down beat point $y_2'$ and discuss whether is a down beat point or not. By finiteness, there exists a down beat point $y_n'$ such that $y_n'<y_{n-1}'<\cdots< y_1'<x$ and $r(y_n')\neq y_n'$.
\end{proof}

This result is also illustrated in Figure \ref{fig:todos_semiflujos}. The semiflows that move the only potential down beat point of $X$, which is $A$, also move $B$, $C$ or both $B$ and $C$ at the same time. Furthermore, this example proves that there are semiflows  that can move a potential down beat point $x\in X$ and fix a down beat point in $U_x$.

\begin{prop}\label{number_height_1} Let $X$ be a finite topological space such that every potential down beat point $x\in X$ has height $1$ and let $n$ be the number of potential down beat points of $X$. Then $ S_F(X)=2^n$.
\end{prop}
\begin{proof}
By hypothesis, $x$ is a potential down beat point if and only if $x$ is a down beat point. For each down beat point $x$, there is a unique semiflow $\varphi$ such that $\varphi(t,x)\neq x$ for every $t>0$ and $\varphi(t,y)=y$ otherwise. By combining these semiflows we can obtain every semiflow defined on $X$. Specifically, we have $2^n$ distinct semiflows. 
\end{proof}

Proposition \ref{number_height_1} shows that not every natural number can be realized as $S_F(X)$ for some finite topological space $X$ of height one. Additionally, in the proof of Theorem \ref{thm_realizacion_todos_los_enteros_positivos}, the finite topological space $X_n$ has height $n+1$. This naturally leads to the following questions: \emph{is there a natural number $m$ such that every positive number can be realized as $S_F(X)$ for some finite topological space $X$ with $ht(X)\leq m$?} \emph{Is there a simple closed formula for $S_F(X)$ in terms of the number of potential down beat points of $X$?}

\bibliography{bibliografia}

\begin{thebibliography}{1}

\bibitem{alexandroff1937diskrete}
P.~S. Alexandroff.
\newblock Diskrete {R}äume.
\newblock {\em Mathematiceskii Sbornik (N.S.)}, 2:501--518, 1937.

\bibitem{barmak2011algebraic}
J.~A. Barmak.
\newblock {\em Algebraic topology of finite topological spaces and
  applications. \textnormal{Lecture Notes in Mathematics 2032}}.
\newblock Springer, 2011.

\bibitem{barmak2020Lefschetz}
J.~A. Barmak, M.~Mrozek, and T.~Wanner.
\newblock A {L}efschetz fixed point theorem for multivalued maps of finite
  spaces.
\newblock {\em Math. Z.}, 294:1477--1497, 2020.

\bibitem{barmak2024conley}
J.~A. Barmak, M.~Mrozek, and T.~Wanner.
\newblock Conley index for multivalued maps on finite topological spaces.
\newblock {\em Found. Comput. Math.}, 2024.

\bibitem{chocano2022dynamical}
P.~J. Chocano, D.~Mond\'ejar~Ruiz, M.~A. Mor\'on, and F.~R. Ruiz~del Portal.
\newblock On the triviality of flows in {A}lexandroff spaces.
\newblock {\em Topol. Appl.}, 339, Part A, 2023.

\bibitem{chocano2020computational}
P.~J. Chocano, M.~A. Mor\'on, and F.~R. Ruiz~del Portal.
\newblock Computational approximations to compact metric spaces.
\newblock {\em Phys. D: Nonlinear Phenom}, 433(133168), 2022.

\bibitem{chocano2024coincidence}
P.~J. Chocano, M.~A. Mor\'on, and F.~R. Ruiz~del Portal.
\newblock Coincidence theorems for finite topological spaces.
\newblock {\em Topol. Methods Nonlinear Anal.}, (To appear), 2024.

\bibitem{mccord1966singular}
M.~C. McCord.
\newblock Singular homology groups and homotopy groups of finite topological
  spaces.
\newblock {\em Duke Math. J.}, 33:465--474, 1966.

\bibitem{stong1966finite}
R.~E. Stong.
\newblock Finite topological spaces.
\newblock {\em Trans. Am. Math. Soc.}, 123(2):325--340, 1966.

\end{thebibliography}
\bibliographystyle{plain}

\newcommand{\Addresses}{{
  \bigskip
  \footnotesize

  \textsc{ P.J. Chocano, Departamento de Matemática Aplicada,
Ciencia e Ingeniería de los Materiales y
Tecnología Electrónica, ESCET
Universidad Rey Juan Carlos, 28933
Móstoles (Madrid), Spain}\par\nopagebreak
  \textit{E-mail address}:\texttt{ pedro.chocano@urjc.es}

}}

\Addresses

\end{document}